\documentclass{amsart}
\usepackage{amssymb,latexsym,amsmath,amsfonts,hhline}
\usepackage{pdfsync}
\input xy
\xyoption{all}

\makeatletter
\newcommand{\leqnomode}{\tagsleft@true}
\newcommand{\reqnomode}{\tagsleft@false}
\makeatother

\newcommand{\cal}{\mathcal}

\makeatletter
\renewcommand{\subsection}{\@startsection{subsection}{2}{0mm}{-2mm}{-2mm}{\bf\normalsize}}

\def\sbsnt#1{\subsection{#1}}
\makeatother

\newtheorem{formula}{}[section]
\newtheorem{definition}[formula]{Definition}
\newtheorem{corollary}[formula]{Corollary}
\newtheorem{remark}[formula]{Remark}
\newtheorem{example}[formula]{Example}
\newtheorem{lemma}[formula]{Lemma}
\newtheorem{theorem}[formula]{Theorem}
\newtheorem{proposition}[formula]{Proposition}

\def\thrm{\begin{theorem}}
\def\thrml#1{\begin{theorem}\label{#1}}
\def\ethrm{\end{theorem}}
\def\prpstn{\begin{proposition}}
\def\prpstnl#1{\begin{proposition}\label{#1}}
\def\eprpstn{\end{proposition}}
\def\rmrk{\begin{remark}}
\def\rmrkl#1{\begin{remark}\label{#1}}
\def\ermrk{\end{remark}}
\def\dfntn{\begin{definition}}
\def\dfntnl#1{\begin{definition}\label{#1}}
\def\edfntn{\end{definition}}
\def\nmrt{\begin{enumerate}}
\def\enmrt{\end{enumerate}}
\def\tm#1{\item[{\rm (#1)}]}
\def\qtn{\begin{equation}}
\def\qtnl#1{\begin{equation}\label{#1}}
\def\eqtn{\end{equation}}
\def\xmpl{\begin{example}}
\def\xmpll#1{\begin{example}\label{#1}}
\def\exmpl{\end{example}}
\def\lmm{\begin{lemma}}
\def\lmml#1{\begin{lemma}\label{#1}}
\def\elmm{\end{lemma}}
\def\crllr{\begin{corollary}}
\def\crllrl#1{\begin{corollary}\label{#1}}
\def\ecrllr{\end{corollary}}
\def\css{\begin{cases}}
\def\ecss{\end{cases}}

\def\proof{\noindent{\bf Proof}.\ }

\def\cD{{\cal D}}

\def\cM{{\cal M}}

\def\cX{{\cal X}}
\def\cY{{\cal Y}}

\def\mZ{{\mathbb Z}}

\DeclareMathOperator{\aut}{Aut}

\DeclareMathOperator{\Fix}{Fix}

\DeclareMathOperator{\id}{id}

\DeclareMathOperator{\iso}{Iso}

\DeclareMathOperator{\orb}{Orb}

\DeclareMathOperator{\sz}{Sz}
\DeclareMathOperator{\syl}{Syl}
\DeclareMathOperator{\sym}{Sym}

\def\eprf{\hfill$\square$}

\def\qaq{\quad\text{and}\quad}
\def\qoq{\quad\text{or}\quad}
\def\ov{\overline}
\newcommand{\grp}[1]{\langle {#1}\rangle}

\def\VRT#1{*=<5mm>[o][F-]{#1}}
 
\def\grphp#1{$\xymatrix@R=10pt@C=10pt@M=0pt@L=2pt{#1}$}

\begin{document}
\title{Coherent configurations associated with TI-subgroups}
\author{Gang Chen}
\address{School of Mathematics and Statistics, Central China Normal University, Wuhan, China}
\email{chengangmath@mail.ccnu.edu.cn}
\thanks{The work of the first author was supported by the 
NSFC (No. 11571129)}
\author{Ilia Ponomarenko}
\address{Steklov Institute of Mathematics at St. Petersburg, Russia}
\email{inp@pdmi.ras.ru}
\date{}

\begin{abstract}
Let $\cX$ be a coherent configuration associated with a transitive 
group~$G$. In terms of the intersection numbers of $\cX$, a necessary
condition for the point stabilizer of $G$ to be a TI-subgroup,
is established. Furthermore, under this condition, $\cX$ is determined 
up to isomorphism by the intersection numbers. It is also proved
that asymptotically, this condition is also sufficient. More precisely, 
an arbitrary homogeneous coherent configuration satisfying this condition 
is  associated with a transitive group, the point stabilizer of which
is a TI-subgroup. As a byproduct of the developed theory, recent results 
on pseudocyclic and quasi-thin association schemes are generalized and
improved. In particular, it is shown that any scheme of prime degree $p$ and 
valency $k$ is associated with a transitive group, whenever
$p>1+6k(k-1)^2$.\\[5mm]
\end{abstract}

\maketitle

\section{Introduction}
Recall that a subgroup $H$ of a finite group $G$ is said to be TI
(trivial intersection) if $H^g\cap H=H$ or $1$ for all $g\in G$. 
In particular, this condition is obviously satisfied if $H$ is of prime 
order or $G$ is
a Frobenius group and $H$ is a Frobenius complement. Some other
examples can be found, e.g., in \cite{BG,HT,W81}. One of the goals of 
this paper is to prove that, at least asymptotically, the concept
of TI-subgroup can be defined in a purely combinatorial way (see
Theorem~\ref{240116g} below). This enables us to generalize and
improve recent results on pseudocyclic and quasi-thin association
schemes~\cite{MP12a,MP12b}.\medskip 

Let $G\le\sym(\Omega)$ be a transitive permutation group and $H$ the 
one point stabilizer of $G$. Suppose that $H$ is a TI-subgroup 
of $G$. Then from Theorem~\ref{200316e} below, it follows that 
\qtnl{260316b}
|X|\in\{1,k\},\qquad X\in\orb(G,\Omega),
\eqtn
where $k$ is the order of $H$.
Moreover, set $c$ to be the maximum number of points fixed by the
permutations belonging to a right $H$-coset of $G$ other than $H$. Then by
Lemma~\ref{190316r}, we have
\qtnl{170116b}
c\le mk,
\eqtn
where  $m$ is the index of $H$ in $N_G(H)$. In 
Section~\ref{260316a}, it is shown that this bound is tight, but in general
conditions~\eqref{260316b} and~\eqref{170116b} do not guarantee
that $H$ is a TI-subgroup of~$G$. Nevertheless, the situation changes if
the ratio $n/m$ with $n=|\Omega|$ is enough large in 
comparison with~$k$ (Proposition~\ref{200316c}).\medskip

What we said in the above paragraph has a clear 
combinatorial meaning. Namely, let $\cX=(\Omega,S)$ be the (homogeneous)
coherent configuration associated with the group~$G$ (the exact
definitions concerning coherent configurations can be found in
Section~\ref{190116a}). We say that $\cX$ is a {\it TI-scheme} if $H$
is a TI-subgroup of~$G$. In this case, condition~\eqref{260316b} means that
\qtnl{260316c}
S=S_1\,\cup\,S_k,
\eqtn
where  $S_1$ and $S_k$ are the sets of basis relations of $\cX$ 
that have valences~$1$ and $k$, respectively.\footnote{In the terminology
of~\cite{Z05}, the set $S_1$ is the thin radical of~$\cX$.} The numbers
$c$, $m$, and $k$ in condition~(2) are, respectively, the 
indistinguishing number\footnote{In the complete colored graph 
representing $\cX$, $c$ is the  maximum number of triangles with fixed 
base, the other two sides of which are 	monochrome arcs (see also~Subsection~\ref{230415b}).}, the cardinality of $S_1$, and the
maximal valency of~$\cX$. Thus, both conditions~\eqref{260316b}
and~\eqref{170116b} can be expressed in terms of~$\cX$ without
mentioning~$G$.

\dfntnl{200316a}
A homogeneous coherent configuration $\cX$ satisfying conditions
\eqref{260316c} and~\eqref{170116b} is said to be a pseudo-TI scheme
of valency~$k$ and index~$n/m$.
\edfntn

The above discussion shows that every TI-scheme
is also pseudo-TI. Furthermore, the result~\cite[Theorem~3.2]{MP12a}
shows that every pseudocyclic scheme is pseudo-TI with $m\in\{1,n\}$ 
and $c=k-1$. On the other hand, many quasi-thin schemes
are also pseudo-TI but this time with $k=2$ (Theorem~\ref{220316a}).
Some other examples of pseudo-TI schemes (in particular, those which 
are not TI) can be found in Section~\ref{260316a}.\medskip

As was mentioned, the main result of the present paper shows that asymptotically,
every pseudo-TI scheme $\cX$ is also TI. However, we can say more. Namely,
in this case, $\cX$ is determined up to isomorphism by its intersection
numbers, i.e., is {\it separable} in terms of~\cite{EP09}. Thus,
at least asymptotically, there is a 1-1 correspondence between the
permutation groups with TI point stabilizers and pseudo-TI schemes.

\thrml{240116g}
Every pseudo-TI scheme of valency $k$ and index greater than $1+6k^2(k-1)$,
is a separable TI-scheme.
\ethrm

The bound in Theorem~\ref{240116g} is not tight. To see this, let
$\cX$ be a pseudo-TI scheme with $k=2$. Then $\cX$ is a quasi-thin
association scheme. Assume that  the index of~$\cX$ is equal to~$2$
(such schemes do exist).  Then from \cite[Theorem~1.1]{MP12b}, 
it follows that $\cX$ is separable and {\it schurian}, i.e., is 
associated with an appropriate permutation group. Thus,
this scheme is TI. However, in this case the bound given in Theorem~\ref{240116g} is not achieved.\medskip

Formula \eqref{260316c} shows that $|S|=m+(n-m)/k$.
Therefore Theorem~\ref{240116g} immediately follows from 
inequality~\eqref{170116b} and the theorem below, in which the numbers
$k$, $c$, and $m$ are called the {\it parameters} of the pseudo-TI scheme~$\cX$.

\thrml{170116c}
Let $\cX=(\Omega,S)$ be a pseudo-TI scheme with parameters $k$, $c$, and~$m$. 
Suppose that 
$$
|S|>m+6c(k-1), 
$$
Then $\cX$ is schurian and separable.
\ethrm

It should be mentioned that according to~\cite[Theorems~1.1 and~1.3]{MP12a},
the statement of Theorem~\ref{170116c} is true if the scheme $\cX$ 
is pseudocyclic and its rank $|S|$ is greater than $ak^4$ for some $a>0$. 
In this case, $\cX$ must be a Frobenius scheme, i.e., $\cX$ is
associated with a Frobenius group,  the point stabilizer of which
coincides with the Frobenius complement. Moreover, since here
$c=k-1$, we can substantially improve the above bound as follows.

\crllrl{170116d}
Any pseudocyclic scheme of valency $k>1$ and rank greater that $1+6(k-1)^2$ is a Frobenius scheme.
\ecrllr

The brilliant Hanaki-Uno theorem on association schemes
 of prime degree~$p$
states that such a scheme must be pseudocyclic~\cite{HU}. In
particular, formula~\eqref{260316c} holds for some $k$ called
the valency of the scheme. Thus, as an immediate consequence
of Theorem~\ref{170116c}, we obtain the following statement.

\crllrl{170116e}
Any scheme of prime degree $p$ and valency $k$ is schurian, whenever $p>1+6k(k-1)^2$.
\ecrllr


To make the paper possibly self-contained, we cite the basics of
coherent configurations in Section~\ref{190116a}. In
Section~\ref{260316a}, we study TI- and pseudo-TI coherent 
configurations in detail; our presentation here includes also
several examples. The proof of Theorem~\ref{170116c} is based on a theory of matching
configurations, which are introduced and studied in Section~\ref{030416a}.
This theory is used in Section~\ref{220216c}, where we prove a quite
general Theorem~\ref{240116a}, which shows that the one point extension of
a homogeneous coherent configuration with sufficiently many relations 
of maximum valency must contain a big matching subconfiguration.
These results are the main ingredients in our proof of 
Theorem~\ref{170116c} given in Section~\ref{090316r}.
\medskip

{\bf Notation.}
Throughout the paper, $\Omega$ denotes a finite set.

The diagonal of the Cartesian product $\Omega\times\Omega$ is denoted by~$1_\Omega$; for any $\alpha\in\Omega$
we set $1_\alpha=1_{\{\alpha\}}$.

For $r\subseteq\Omega\times\Omega$, set $r^*=\{(\beta,\alpha):\ (\alpha,\beta)\in r\}$ and
$\alpha r=\{\beta\in\Omega:\ (\alpha,\beta)\in r\}$ for all $\alpha\in\Omega$.

For $s\subseteq \Omega\times\Omega$ we set
$r\cdot s=\{(\alpha,\gamma):\ (\alpha,\beta)\in r,\ (\beta,\gamma)\in s$
for some $\beta\in\Omega\}$. If $S$ and $T$ are sets of relations, we set
$S\cdot T=\{s\cdot t:\ s\in S,\, t\in T\}$.

For $S\in 2^{\Omega^2}$, set $S^\#=S\setminus\{\varnothing\}$,
denote by $S^\cup$ the set of all unions of the elements of $S$, 
and put $S^*=\{s^*:\ s\in S\}$
and $\alpha S=\cup_{s\in S}\alpha s$, where $\alpha\in\Omega$.


\section{Coherent configurations}\label{190116a}

This section contains well-known basic facts on coherent configurations.
All of them can be found in~\cite{EP09} and papers cited there.

\sbsnt{Main definitions.}\label{230415b}
Let $\Omega$ be a finite set, and let $S$ be a partition of a subset
of~$\Omega\times\Omega$. Then the pair $\cX=(\Omega,S)$ is called
a {\it partial coherent configuration} on $\Omega$ if the following 
conditions hold:
\nmrt
\tm{C1} $1_\Omega\in S^\cup$, 
\tm{C2} $S^*=S$, 
\tm{C3} given $r,s,t\in S$, the number $c_{rs}^t=|\alpha r\cap\beta s^*|$
does not depend on the choice of $(\alpha,\beta)\in t$. 
\enmrt
If $S$ is a partition of $\Omega\times\Omega$, then $\cX$ is a {\it coherent
configuration} in the usual sense.
The elements of $\Omega$ and $S$, and the numbers $c_{rs}^t$ are
called the {\it points} and {\it basis relations}, and the {\it intersection numbers} of~$\cX$, respectively. The numbers $|\Omega|$ 
and $|S|$ are called the {\it degree} and the {\it rank} of~$\cX$. 
The basis relation containing a pair
$(\alpha,\beta)\in\Omega\times\Omega$ is denoted by $r(\alpha,\beta)$.
For $r,s\in S^\cup$, the set of basic relations contained in $r\cdot s$  
is denoted by~$rsv$. 
\medskip

Denote by $\Phi=\Phi(\cX)$ the set of $\Lambda\subseteq\Omega$ such that
$1_\Lambda\in S$. The elements of $\Phi$ are called the {\it fibers} 
of~$\cX$. In view of condition (C1), the set $\Omega$ is the disjoint 
union of all of them. Moreover, for each $r\in S$, there exist uniquely
determined fibers $\Lambda$ and $\Delta$ such that
$r\subseteq\Lambda\times\Delta$. Note that
the number $|\alpha r|=c_{rr^*}^t$ with $t=1_\Lambda$,
does not depend on $\alpha\in\Lambda$. This number is called the 
{\it valency} of~$r$ and denoted~$n_r$. One can see that 
if $\Delta\in\Phi^\cup$ and $S_\Delta$ is the set of all basis
relations contained in $\Delta\times\Delta$, then 
$$
\cX_\Delta=(\Delta,S_\Delta)
$$ 
is a partian coherent configuration, called the {\it restriction}
of~$\cX$ to~$\Delta$.

\sbsnt{Homogeneous coherent configurations.}
A coherent configuration $\cX$ is said to be {\it homogeneous} or a 
{\it scheme} if $1_\Omega\in S$. In this case,
$n_r=n_{r^*}=|\alpha r|$ for all $r\in S$ and $\alpha\in\Omega$. Given a positive integer $k$, we set
\qtnl{240116c}
S_k=\{r\in S:\ n_r=k\}.
\eqtn
One can see that the set $S_1$ forms a group with respect to the relation
product. In what follows, we will also use the following well-known
identities for the intersection numbers of a homogeneous coherent
configurations:
\qtnl{150410a}
n_rn_s=\sum_{t\in S}n_tc_{rs}^t,\qquad
c_{r^*s^*}^{t^*}=c_{sr}^t,\qquad
n_tc_{rs}^{t^*}=n_rc_{st}^{r^*}=n_sc_{tr}^{s^*}
\eqtn
for all $r,s,t\in S$.\medskip

The {\it indistinguishing number} $c(r)$ of a relation $r\in S$ is defined
to be the sum of the intersection numbers $c_{ss^*}^r$ 
with $s\in S$ (see \cite{MP12a}). For each pair $(\alpha,\beta)\in r$, 
we obviously have
\qtnl{100814c}
c(r)=|\Omega_{\alpha,\beta}|,\quad\text{where}\quad \Omega_{\alpha,\beta}=\{\gamma\in\Omega:\ r(\gamma,\alpha)=r(\gamma,\beta)\}.
\eqtn
The maximum $c(\cX)$ of the numbers $c(r)$, $r\in S^\#$, is called
the {\it indistinguishing number} of the coherent configuration~$\cX$.


\sbsnt{Point extensions.}\label{300116a}
There is a natural partial order\, $\le$\, on the set of all coherent configurations on the same set.
Namely, given two coherent configurations $\cX=(\Omega,S)$ and
$\cX'=(\Omega,S')$, we set
$$
\cX\le\cX'\ \Leftrightarrow\ S^\cup\subseteq (S')^\cup.
$$
The minimal and maximal elements with respect to this ordering are,
respectively, the {\it trivial} and {\it complete} coherent
configurations: the basis relations of the former one are the reflexive
relation $1_\Omega$ and (if $n>1$) its complement in $\Omega\times\Omega$,
whereas the basis relations of the latter one are all singletons.\medskip

Given two coherent configurations $\cX_1$ and $\cX_2$ on $\Omega$, there is a uniquely determined coherent configuration
$\cX_1\cap\cX_2$ also on $\Omega$, the relation set of which is $(S_1)^\cup\cap(S_2)^\cup$, where $S_i$ is the set of basis
relations of~$\cX_i$, $i=1,2$. This enables us to define the {\it coherent closure}
of a set $R$ of binary relations on $\Omega$ by
$$
(\Omega,\ov R)=\bigcap_{T:\ R\subseteq T^\cup}(\Omega,T),
$$
where $T$ runs over the partitions of $\Omega\times\Omega$ such that 
$(\Omega,T)$ is a coherent configuration. Note that if $(\Omega,R)$ is
a coherent configuration, then $\ov R=R$.\medskip

Let $\cX=(\Omega,S)$ be a coherent configuration. For the points
$\alpha,\beta,\ldots$, denote by $\cX_{\alpha,\beta,\ldots}$ the
coherent closure of the set $R=S\,\cup\,\{1_\alpha,1_\beta,\ldots\}$.
This coherent configuration is called the
{\it point extension}  of $\cX$ with respect to the points
$\alpha,\beta,\ldots$.  In other words,
$\cX_{\alpha,\beta,\ldots}$ is the smallest coherent configuration 
on $\Omega$ that is larger than or equal to $\cX$ and has the singletons
$\{\alpha\},\{\beta\},\ldots$ as fibers. The following lemma
immediately follows from the definitions.

\lmml{260511a}
Let $\Phi_\alpha$ and $S_\alpha$ be, respectively, the sets of fibers and
basis relations of the coherent configuration $\cX_\alpha$, where
$\alpha\in\Omega$. Then 
$$
\alpha r\in(\Phi_\alpha)^\cup\qaq r_{x,y}\in (S_\alpha)^\cup,\qquad r,x,y\in S,
$$
where $r_{x,y}=r\,\cap(\alpha x\times\alpha y)$.  Moreover, 
$|\beta r_{x,y}|=c_{r^{}y^*}^{x^*}$ for all $\beta\in\alpha x$.
\elmm

In what follows, we often use the notation $r_{x,y}$ defined in 
Lemma~\ref{260511a}. Note here that $(r_{x,y})^*=(r^*)_{y,x})$
for all $r,x,y$.

\sbsnt{Direct sum and tensor product.}
Let $\cX=(\Omega,S)$ and $\cX'=(\Omega',S')$ be coherent configurations. Denote by
$\Omega\sqcup\Omega'$ the disjoint union of~$\Omega$ and~$\Omega'$, and by
$S\boxplus S'$ the union of the set $S\sqcup S'$ and the set of all relations
$\Delta\times\Delta'$ and $\Delta'\times\Delta$ with $\Delta\in\Phi(\cX)$ 
and $\Delta'\in\Phi(\cX')$. Then the pair
$$
\cX\boxplus\cX'=(\Omega\sqcup\Omega',S\boxplus S')
$$
is a coherent configuration called the {\it direct sum} of~$\cX$ and~$\cX'$.
One can see that $\cX\boxplus\cX'$ equals the coherent closure of 
the set $S\sqcup S'$.
\medskip

Set $S\otimes S'=\{u\otimes u':\ u\in S,\ u'\in S'\}$, where $u\otimes u'$ is
the relation on $\Omega\times\Omega'$ consisting of all 
pairs $((\alpha,\alpha'),(\beta,\beta'))$
with $(\alpha,\beta)\in u$ and $(\alpha',\beta')\in u'$. Then the pair
$$
\cX\otimes\cX'=(\Omega\times\Omega',S\otimes S')
$$
is a coherent configuration called the {\it tensor product} of~$\cX$ and~$\cX'$.
One can see that $\cX\otimes\cX'$ equals the coherent closure of the set
$S\otimes 1_{\Omega'}\,\cup\, 1_\Omega\otimes S'$.

\sbsnt{Isomorphisms and schurity.}
Two partial coherent configurations are called {\it isomorphic} if 
there exists a bijection between their point sets that induces a bijection
between their sets of basis relations. Each such bijection is called an 
{\it isomorphism} between 
these two configurations. The group of all isomorphisms of a partial coherent configuration $\cX=(\Omega,S)$ to itself 
contains a normal subgroup
$$
\aut(\cX)=\{f\in\sym(\Omega):\ s^f=s,\ s\in S\}
$$
called the {\it automorphism group} of~$\cX$, where
$s^f=\{(\alpha^f,\beta^f):\ (\alpha,\beta)\in s\}$.\medskip

Conversely, let $G\le\sym(\Omega)$ be a permutation group, and let $S$ 
be a set of orbits of the component-wise action of $G$ 
on~$\Omega\times\Omega$. Assume that $S$ satisfies conditions~(C1) and~(C2).
Then, $\cX=(\Omega,S)$ is a partial coherent configuration;
we say that $\cX$ {\it is  associated} with~$G$. A partial coherent configuration on $\Omega$ is said to be {\it schurian} if it is 
associated with some permutation group on~$\Omega$. It is easily seen that a partial coherent configuration~$\cX$ is schurian 
if and only if it is the partial coherent configuration associated
with the group~$\aut(\cX)$.

\sbsnt{Algebraic isomorphisms and separability.}
Let $\cX=(\Omega,S)$ and $\cX'=(\Omega',S')$ be partial coherent configurations. A bijection $\varphi:S\to S',\ r\mapsto r'$ is 
an {\it algebraic isomorphism} from~$\cX$ onto~$\cX'$ if
\qtnl{f041103p1}
c_{rs}^t=c_{r's'}^{t'},\qquad r,s,t\in S.
\eqtn
In this case, $\cX$ and $\cX'$ are said to be {\it algebraically isomorphic}. Each isomorphism~$f$ from~$\cX$ onto~$\cX'$ 
induces an algebraic isomorphism $\varphi_f:r\mapsto r^f$
between these configurations. The set of all
isomorphisms inducing the algebraic isomorphism~$\varphi$ is denoted by $\iso(\cX,\cX',\varphi)$. In particular,
\qtnl{190316b}
\iso(\cX,\cX,\id_S)=\aut(\cX),
\eqtn
where $\id_S$ is the identity mapping on $S$. A partial coherent configuration~$\cX$ is said to be {\it separable} if for any algebraic 
isomorphism~$\varphi:\cX\to\cX'$, the set $\iso(\cX,\cX',\varphi)$ is not empty.\medskip

An algebraic isomorphism $\varphi$ induces a bijection from $S^\cup$
onto $(S')^\cup$: the union $r\cup s\cup\cdots$ of basis relations 
of $\cX$ is taken to $r'\cup s'\cup\cdots$. This bijection is also
denoted by $\varphi$. One can see that $\varphi$ preserves the reflexive
basis relations. This extends $\varphi$ to a bijection
$\Phi(\cX)\to \Phi(\cX')$ so that $(1_{\Delta^{}})'=1_{\Delta'}$.\medskip

In the above notation, assume that $\cY\ge\cX$ and $\cY'\ge\cX'$ 
are partial coherent configurations. We say that the algebraic isomorphism
$\psi:\cY\to\cY'$ {\it extends} $\varphi$ if $r^\psi=r^\varphi$ 
for all $r\in S$. In this case, we obviously have
$$
\iso(\cY,\cY',\psi)\subseteq\iso(\cX,\cX',\varphi).
$$
This immediately proves the following statement (take $\cY=\cX_\alpha$).

\lmml{150316d}
Let $\cX$ be a coherent configuration and $\alpha\in\Omega$.
Assume that for every algebraic isomorphism $\varphi:\cX\to\cX'$,
there exists  an algebraic isomorphism
$\varphi_{\alpha,\alpha'}:\cX^{}_{\alpha^{}}\to\cX'_{\alpha'}$ 
with $\alpha'\in\Omega'$ that extends $\varphi$. Then $\cX$ is
separable if  so is $\cX_\alpha$.
\elmm

\sbsnt{1-regular coherent configurations.}
A point $\alpha\in\Omega$ of the coherent configuration $\cX$ is called {\it regular} if
$$
|\alpha r|\le 1\quad\text{for all}\ \,r\in S. 
$$
Obviously, every point of a fiber $\Delta\in\Phi(X)$ is regular, whenever
$\Delta$ contains at least one regular point. Therefore,
the set of all regular points of $\cX$ is the union of fibers. If this 
set is not empty, then the coherent configuration 
$\cX$ is said to be {\it 1-regular}; we say that $\cX$ is 
{\it semiregular} if each of its points is regular.
Note that in the homogeneous case, a coherent configuration is 1-regular if and 
only if it is a thin  scheme in the sense of~\cite{Z05}. The 
following statement can be found in~\cite[Theorem~9.3]{EP09}.

\thrml{150316c}
Every 1-regular coherent configuration is schurian and separable.\eprf
\ethrm

It is easily seen that a coherent configuration is semiregular if and only if
every basis relation of it is  a {\it matching}, i.e., a binary relation
of the form 
$$
r=\{(\alpha,f(\alpha)):\ \alpha\in\Delta\},\quad\text{where}\ \ f:\Delta\to\Delta'
\ \text{is a bijection.}
$$ 
Clearly, $|r|=|\Delta|=|\Delta'|$. Note that if $r'$ is a matching
with respect to a bijection $f':\Delta'\to\Delta''$, then $r\cdot r'$ is 
a matching with respect to the composition $f\circ f'$.

\section{TI- and pseudo-TI schemes}\label{260316a}

\sbsnt{Parameters.}
Throughout this section, $G\le\sym(\Omega)$ is a transitive group
and $H=G_\alpha$ is the stabilizer of the point $\alpha$ in~$G$.
The following statement gives a necessary and sufficient condition 
for $H$ to be a TI-subgroup of~$G$. 

\thrml{200316e}
Let $\cX=(\Omega,S)$ be the coherent configuration associated 
with the group~$G$, and let $k$ be the maximum of $n_s$, $s\in S$. 
Then $H$ is a TI-subgroup of~$G$ if and only if $S=S_1\cup S_k$ and 
$H$ acts semiregularly on $\alpha S_k$. In particular, $k=|H|$,
whenever $H$ is a TI-subgroup of~$G$.
\ethrm
\proof By the transitivity of~$G$, for every $\beta\in \Omega$ there 
exists $g\in G$ such that $\beta=\alpha^g$. It follows that
\qtnl{200316f}
H_\beta=G_{\alpha,\beta}=G_\alpha\cap G_\beta=H\cap H^g.
\eqtn
Moreover, it is easily seen that that $\beta\in\alpha S_1$ if and
only if $G_\alpha=G_\beta=H^g$. Thus,
\qtnl{200316g}
\beta\in\alpha S_1\qquad\Longleftrightarrow\qquad H=H^g.
\eqtn
Assume that $H$ is a TI-subgroup of~$G$. Then from~\eqref{200316f} 
and~\eqref{200316g}, it follows that $H_\beta=1$ for all 
$\beta$ belonging to $\Omega'=\Omega\setminus\alpha S_1$. Therefore, 
$H$ acts semiregularly on $\Omega'$. Since 
$$
\orb(H,\Omega')=\{\alpha s:\ s\in S\setminus S_1\},
$$
this implies that $n_s=|\alpha s|=|H|$ for all $s\in S\setminus S_1$.
It follows that $|H|=k$ and hence
$S\setminus S_1=S_k$ and $\Omega'=\alpha S_k$. This proves
the ``only if" part. Conversely, assume that $S=S_1\cup S_k$ and 
$H$ acts semiregularly on $\alpha S_k$. Let $g\in G$ be
such that $H\ne H^g$. Then from~\eqref{200316g} with $\beta=\alpha^g$,
it follows that $\beta\in \alpha S_k$. Therefore, $H_\beta=1$ and hence
 $H\cap H^g=1$ in view of~\eqref{200316f}. Thus, $H$ is
a TI-group, as required.\eprf\medskip

From now on, we assume that $H$ is a TI-subgroup of $G$, i.e.,
the coherent configuration $\cX$ is a TI-scheme. Note that the 
group $G$ is not uniquely determined by $\cX$ as the following 
example shows (it also shows that, in general, $G\ne\aut(\cX)$).

\xmpll{210316s}
{\rm Let $C=\mZ_{p^2}$, $H\le\aut(C)$, and $G$ the permutation
group induced by the action of $C\rtimes H$ on $\Omega=C$, in which
$C$ acts by right multiplications and $H$ acts naturally. If now
$|H|=p$, then the coherent configuration $\cX$ is a TI-scheme and
$\aut(\cX)\cong \mZ_p\wr\mZ_p$.
One can see that $G$ coincides with its normalizer in $\aut(\cX)$.
Therefore, there are $p^{p-2}$ conjugates of~$G$ in $\aut(\cX)$. 
The coherent configurations associated with them are equal to~$\cX$.} 
\exmpl

However, from Theorem~\ref{200316e}, it follows that the orders of 
the groups $G$ and $H$ are uniquely determined by the coherent
configuration~$\cX$, namely, $|G|=nk$ and $|H|=k$, where $n=|\Omega|$. 
Of course, $\cX$ also determines the rank $r=|S|$ of the permutation
group~$G$ and the order of the group $N=N_G(H)$, where the latter
immediately follows from~\eqref{200316g} showing that $|N:H|=m$ with
$m=|S_1|$. Note that 
$$
n=m+(r-m)k,
$$
because $S=S_1\,\cup\,S_k$ by Theorem~\ref{200316e}.
Nevertheless, not every schurian homogeneous coherent configuration 
with two distinct valences of basis relations is a TI-scheme.

\xmpl
{\rm Let $\cX$ be the tensor product of a regular scheme on $m$ points
and the Johnson scheme $J=J(7,2)$, which is a coherent configuration 
of rank~$3$ on~$21$ points. Then the valences of $\cX$ are $1$ and $10$,
and
$$
n=21m,\quad k=10,\quad r=3m.
$$ 
It is easily seen that $\cX$ is a TI-scheme if and only if so is $J$.
However, if $J$ is a TI-scheme, then by Theorem~\ref{200316e}, the
group $\aut(J)$ contains a subgroup of order 
$21\cdot 10$ acting on $21$ points as a rank~$3$ group. But such a group 
must be solvable and primitive. Thus, the number~$21$ 
is a prime power, contradiction.} 
\exmpl

In fact, this example is a slight generalization of an example arising 
in studying pseudocyclic schemes \cite[p.48]{BCN}. As in that case,
it is quite natural to introduce into consideration the indistinguishing
number defined in~\eqref{100814c}. 

\lmml{190316r}
In the above notation, let $c$ be the indistinguishing number of
the TI-scheme~$\cX$. Then $c\le mk$ and this bound is tight.
\elmm
\proof Denote by $\Fix(x)$ the set of all points fixed by a permutation
$x\in G$. Clearly, $\Fix(1)=\Omega$. We claim that
\qtnl{070416a}
|\Fix(x)|\in\{0,m\}
\eqtn
for all non-identity permutations $x\in G$. Indeed, $\Fix(x)=\varnothing$ if
$x$ belongs to no $G_\beta$, $\beta\in\Omega$. Suppose that the set 
$\Fix(x)$ is not empty. Then $x\in G_\beta$ for some~$\beta$. By the
transitivity of $G$, without loss of generality, we may 
assume that $\beta=\alpha$. Then
$$
\alpha S_1\subseteq \Fix(x),
$$ 
see formula~\eqref{200316g}.
Moreover, by Theorem~\ref{200316e}, the group $H=G_\alpha$ acts
semiregularly on the set $\Omega\setminus \alpha S_1$. This implies that
$|\Fix(x)|=|\alpha S_1|=m$,
which proves claim~\eqref{070416a}.\medskip

On the other hand, in view 
of \cite[Lemma~2.3]{PV16}, we have 
\qtnl{100814d}
c=\max_{g\in G\setminus H}|\bigcup_{x\in Hg}\Fix(x)|.
\eqtn
By formula~\eqref{070416a}, this yields
$$
c\le \max_{g\in G\setminus H}\sum_{x\in Hg}|\Fix(x)|\le mk,
$$
as required. Furthermore, the bound is attained if and only if 
there exists a coset $Hg\ne H$ such that 
$$
\Fix(x)\ne\varnothing\qaq\Fix(x)\cap\Fix(y)=\varnothing
$$ 
for all distinct $x,y\in Hg$. Let us identify the points of $\Omega$ with 
the right $H$-cosets of~$G$ so that $G$ acts on $\Omega$ by the right multiplications.
Then the first condition holds if $Hzx=Hz$ for some $z\in G$, or, equivalently,
if $x^G$ intersects~$H$. Next, the second condition is always true, 
whenever $|H|=2$: indeed, if $Hz\in \Fix(x)\cap\Fix(y)$ for some $z\in G$,
then 
$$
Hzx=Hz=Hzy,
$$
which implies that the point stabilizer of $Hz$
contains $1$, $x$, and $y$, and hence $x=y$. Thus, the bound is attained 
if, for instance, 
$$
G=\sym(6),\qquad H=\grp{(1,2)(3,4)},\qquad g=(1,3)(2,4).
$$
In this case, the coset $Hg$ consists of $x=g=(1,3)(2,4)$ and
$y=(1,4)(2,3)$, and hence $x^z\in H$ for $z=(2,3)$ and
$y^z\in H$ for $z=(2,4)$.\eprf\medskip

From Theorem~\ref{200316e} and Lemma~\ref{190316r}, we immediately 
get the following statement, which, in some sense, justifies 
Definition~\ref{200316a}.

\thrml{200316b}
Every TI-scheme is pseudo-TI.
\ethrm

In general, the converse statement is not true, 
see Example~\ref{230316a} below. Theorem~\ref{240116g}, which
is proved later, shows that it becomes true if the index
$n/m$ is sufficiently large in comparison with~$k$. The
bound established in that theorem can be slightly improved
if the pseudo-TI scheme in question is assumed to be schurian. 

\prpstnl{200316c}
Every schurian pseudo-TI scheme of index greater than $2k(k-1)$
is TI.
\eprpstn
\proof Let $\cX$ be a pseudo-TI scheme. Then $k$ is the maximum valency
of~$\cX$, $c\le mk$, and, in view of the hypothesis, $n/m>2k(k-1)$.
It follows that 
$$
n>2mk(k-1)\ge 2c(k-1).
$$
Assume that $\cX$ is associated with a group $G$. Then
according to \cite[Theorem~3.1 and formula~(2)]{PV16}, this implies
that $H_\beta=G_{\alpha,\beta}=1$ for some point $\beta$.
Note that  $\beta\not\in\alpha S_1$, for otherwise $H=1$.
Consequently, 
$$
|H|=|\alpha s|=k,
$$ 
where $s=r(\alpha,\beta)$. Since $S=S_1\cup S_k$, 
the group $H$ acts semiregularly on the set~$\alpha S_k$. Thus, $H$ is a TI-subgroup
of~$G$ by Theorem~\ref{200316e}, as required.\eprf

\sbsnt{Elementary coset schemes.}\label{240116r}
One of the aims of this subsection is to determine the quasi-thin
schemes, which are pseudo-TI. To this end, let $\cX=(\Omega,S)$ be 
a homogeneous coherent configuration. Assume that there exists a 
set $T\subseteq S_1$ such that
\qtnl{210316r}
ss^*=T\qquad\text{for all}\ \, s\in S\setminus S_1
\eqtn
(see~\cite[Section~6.7]{Z05}). In this case, we say that $\cX$
is an {\it elementary coset scheme}. Special classes of these schemes
have been studied in \cite[Section~3]{M09} and~\cite{EP16}.	 In
both cases, $\cX$ was a wedge product of two regular schemes: one
of them is the restriction of $\cX$ to $\alpha S_1$, and the other
one is  the quotient of $\cX$ modulo~$T$.

\xmpll{210316g}
{\rm Let $C=\mZ_n$ and $H\le\aut(C)$ the group of prime 
order~$p$ dividing~$n$. Denote by $\cX$ the TI-scheme associated 
with the  group $G=C\rtimes H$ (see Example~\ref{210316s}). Thus,
$\Omega=C$ and $S$ consists of the relations
$s_c=(e,c)^G$, where $e$ is the identity of~$C$
and $c\in C$; in particular, 
$$
S_1=\{s_c:\ c\in C\ \text{and}\ |c|\ \,\text{divides}\ n/p\},
$$ 
where $|c|$ is the order of $c$. It is straightforward to check that 
$\cX$ is an elementary coset  
scheme with $T=\{s_c:\ |c|=p\}$, see~\cite[Subsection~8.2]{EP16}.}
\exmpl

Let $\cX$ be an elementary coset scheme.
From \cite[Lemma~6.7.1]{Z05}, it follows that $n_s=n_{ss^*}$ for all 
$s\in S$. Since $T\subseteq S_1$, this implies that $S=S_1\cup S_k$, where
$k=|T|$.  Moreover, in view of~\eqref{210316r}, one can easily find that
$$
c_{ss^*}^r=\css
k    &\text{if $r\in S_1$ and $s\in S_k$,}\\
0    &\text{if $r\in S_k$,}\\
\delta_{r,1}	&\text{if $r,s\in S_1$,}\\
\ecss
$$
where $\delta_{r,1}$ is the Kronecker symbol. It immediately follows that 
$c(r)=0$ for all $r\in S\setminus S_1$ and $c(r)=n-m$ for $r\in (S_1)^\#$,
where $m=|S_1|$. Thus, $c=n-m$.  By the definition of pseudo-TI scheme, 
this proves the following statement.

\prpstnl{210316y}
An elementary coset scheme is pseudo-TI only if $n\le m(k+1)$.
\eprpstn

Note that for the elementary coset scheme from Example~\ref{210316g}, we have 
$k=p$ and  $m=n/p$. In general, we do not know whether the inequality in
Proposition~\ref{210316y} is sufficient for an elementary coset scheme to
be pseudo-TI. However, the following statement shows that this is ``almost" 
true for $k=2$.

\thrml{220316a}
A schurian quasi-thin scheme is not TI if and only if it is elementary
coset scheme and $n\ge 3m$ with possible exception $n=3m$.\footnote{As we
will see below, for
$n=3m$, there exist TI and non-TI quasi-thin elementary coset schemes.}
\ethrm
\proof Let $\cX$ be a schurian quasi-thin scheme; in particular,
$k=2$ and
$$
S=S_1\cup S_2\qaq \orb(G_\alpha)=\{\alpha s:\ s\in S\},
$$
where $G=\aut(\cX)$. The ``if part"
immediately follows from Proposition~\ref{210316y} and Theorem~\ref{200316b}. 
To prove the ``only if" part, assume that $n<3m$ and verify that $\cX$ is a TI-scheme.
Suppose first that $\cX$ has at least two orthogonals,
i.e., the relations $s\in S^\#$ with $c(s)\ne 0$. Then by
\cite[Corollary~6.4]{MP12b}, the coherent configuration
$\cX_\alpha$ is 1-regular. Therefore, the group $G_\alpha$ 
acts semiregularly on $\alpha S_2$. By Theorem~\ref{200316e}, this
implies that $G_\alpha$ is a TI-subgroup of $G$, and hence $\cX$ is
a TI-scheme.\medskip

Now, we may assume that $\cX$ has exactly one
orthogonal, say $u$. If $n_u=2$, then $u$ is the disjoint union of $m$
cliques of size~$3$, and
$$
\cX\cong \cX_1\otimes\cX_2,
$$
where $\cX_1$ and $\cX_2$ are, respectively, trivial and regular 
coherent configurations of degrees~$3$ and~$m$ (see the proof of
Theorem~5.2 in~\cite{MP12b}). It follows that the group
$\aut(\cX)\cong\aut(\cX_1)\times\aut(\cX_2)$ is of order $6m=2n$.
Thus, $|G_\alpha|=|G|/n=2$ and hence $\cX$ is a TI-scheme.\medskip

 Finally, let $n_u=1$. Then condition~\eqref{210316r} is
satisfied for $T=\{1,u\}$. Therefore $\cX$ is an elementary coset
scheme of index~$2$ (the index can not be equal to~$1$, 
because $u\in S_1$). Then one can easily prove that 
$$
\cX_\alpha\cong \cX_1\boxplus \cX_2,
$$ 
where $\cX_1$ is the complete coherent configuration on $\alpha S_1$
and $\cX_2$ is a semiregular coherent configuration on $\alpha S_2$.
This implies that $\cX_\alpha$ is 1-regular and we are done
by Theorem~\ref{200316e} as above.\eprf\medskip

Example~\ref{210316g} for $n=3\cdot 2^n$ and $p=2$, gives a quasi-thin TI-scheme 
with $n=3m$. On the other hand, the following example shows that there are
quasi-thin pseudo-TI-schemes with $n=3m$, which is not TI.

\xmpll{230316a}
{\rm Let $G$ be a unique group of order $2^5\cdot 3^2$ with elementary abelian
socle of order~$8$ (in the GAP notation \cite{gap}, this is the group 
$[288,859]$ with structure description {\rm ``A4\ x\ SL(2,3)"}).  It has
the center of order~$2$ and exactly three non-normal subgroups of order~$4$ that
lie in the socle and do not contain the center. Take one of them, say $H$; the
other two subgroups are conjugate. Then the action of $G$ on the right $H$-cosets
by the right multiplication is a transitive group of degree $n=72$. The coherent
configuration associated with $G$ (in this action) is a quasi-thin elementary
coset scheme $\cX$ with $n=3m$. Therefore, it is pseudo-TI. On the other
hand, $\aut(\cX)=G$ and $G$ contains no subgroup of index~$2$. This implies
that $\cX$ can not be a TI-scheme.}
\exmpl

As it was shown in~\cite{MP12b}, any non-schurian quasi-thin scheme $\cX$
has $2$ or $3$ orthogonals, they generate the Klein group, and 
the index of $\cX$ equals~$4$ or~$7$. It is very likely that such a
scheme is pseudo-TI.

\sbsnt{Further examples of pseudo-TI $p$-schemes.}
In general, it is not easy to find or estimate the indistinguishing
number of a coherent configuration $\cX$ and hence to check whether 
or not $\cX$ is a pseudo-TI scheme. However, if $mk\ge n$, then, obviously, 
inequality~\eqref{170116b} holds. Therefore, $\cX$ is pseudo-TI if and only if $S=S_1\cup S_k$. Such schemes
do exist as the following direct interpretation of a 
construction in~\cite{K15} shows.

\thrml{250316a}
For every prime $p\equiv 3\ (\text{\rm mod}\hspace{2pt}4)$
and each of two non-abelian group~$H$ of order~$p^3$, there
exists a non-schurian pseudo-TI scheme with 
$$
(n,m,k)=(p^3,p^2,p)
$$
admitting $H$ as a regular automorphism group.
\ethrm

Another example comes from the Suzuki group $G=\sz(q)$, where
$q=2^{2n+1}$. The following facts are well known, see, e.g.,
\cite{BG}. The group $G$ is of order $q^2(q^2+1)(q-1)$, has
$q^2+1$ Sylow $2$-subgroups, and any two distinct Sylow $2$-subgroups
intersect at identity. In particular, each $H\in\syl_2(G)$ is 
a TI-subgroup of~$G$. In addition, $|N_G(H):H|=q-1$. Thus, the
coherent configuration associated with the action of $G$ on the
right $H$-cosets by the right multiplication is a TI-scheme with
$$
(n,m,k)=(q^2(q-1),q-1,q^2).
$$
It should be mentioned that $H$ is a non-abelian $2$-group.
Some other examples of simple groups with TI-subgroups can be found
in~\cite{W81}.

\section{Matching configurations}\label{030416a}
Let $\cM=(\Delta,M)$ be a partial coherent configuration, the 
fibers $\Delta_x$ of which are indexed by the elements of a set $X$. Denote
by $M(x,y)$ the set of all its basis relations contained in
$\Delta_x\times\Delta_y$. In what follows, for all $x,y\in X$,
we set $1_x=1_{\Delta_x}$ and write $x\sim y$ if $M(x,y)$ is
a partition of $\Delta_x\times\Delta_y$ into matchings.

\dfntnl{290216b}
We say that $\cM$ is a matching configuration if for every
$x,y\in X$ with $M(x,y)\ne\varnothing$, either $x\sim y$ 
or $x=y$ and $M(x,y)=\{1_x\}$.
\edfntn

Any semiregular coherent configuration is, obviously,
a  matching configuration. In general, with a matching 
configuration $\cM$, we associate an undirected graph $D=D(\cM)$ with 
vertex set $X$, in which the vertices $x$ and $y$ are 
adjacent if and only if $x\sim y$. In particular, we do 
not exclude loops. The following two statements 
describe the adjacency and triangles in $D$:
\qtnl{220216b}
x\sim y\quad\Leftrightarrow\quad
|\Delta_x|=|M(x,y)|=|M(y,x)|=|\Delta_y|,
\eqtn
\qtnl{230116a}
x\sim y\sim z\sim x\quad\Rightarrow\quad
M(x,y)\cdot M(y,z)=M(x,z).
\eqtn
The first of them is obvious, whereas the second one follows from
condition~(C3) in the definition of partial coherent configuration.
\medskip

The matching configuration $\cM$ is said to be {\it saturated} 
if for any set $Y\subseteq X$ with at most $4$ elements, there exists
a vertex of $D$ adjacent with every 
vertex of~$Y$. Note that in this case, any two vertices $x,y\in X$ 
are connected by a $2$-path $P=(x,z,y)$: take $Y=\{x,y\}$.
In particular, $D$ is a graph of diameter at most~$2$.

\thrml{200116a}
Let $\cM=(\Delta,M)$ be a saturated matching configuration, and let
$\cY=(\Delta,T)$ be the coherent closure of $M$. Then 
$T=(M\cdot M)^\#$.
\ethrm
\proof We recall that the fibers $\Delta_x\in\Phi(\cM)$ are indexed 
by the elements~$x\in X$. Since $\cM$ is saturated, 
formula~\eqref{220216b} implies that the number $|\Delta_x|$ does not
depend on $x\in X$; denote this number by~$k$. For any $d$-path
$P=(x_1,\ldots,x_{d+1})$ of the graph $D$, set
$$
M_P=M(x_1,x_2)\cdot\,\cdots\,\cdot M(x_d,x_{d+1}).
$$
Clearly, $M_P$ consists of matchings contained in 
$\Delta_{x_1}\times \Delta_{x_{d+1}}$. 

\lmml{191215b}
For every path $P$, the set $M_P$ is a partition with $k$ classes. Moreover, 
this partition does not depend on the choice of the path $P$ connecting 
$x_1$ and $x_{d+1}$.
\elmm
\proof To prove the first statement, it suffices to verify that  the equalities 
\qtnl{191215e}
a\cdot M(x_2,x_3)=M(x_1,x_2)\cdot M(x_2,x_3)=M(x_1,x_2)\cdot b
\eqtn
hold for any $a\in M(x_1,x_2)$ and $b\in M(x_2,x_3)$. Note that if the first 
equality is true for all $2$-paths $P=(x_1,x_2,x_3)$ and $a$, then 
applying it to the path $P^*=(x_3,x_2,x_1)$, we have
$$
M(x_1,x_2)\cdot b=(b^*\cdot M(x_2,x_1))^*=(M(x_3,x_2)\cdot M(x_2,x_1))^*=M(x_1,x_2)\cdot M(x_2,x_3).
$$
Thus, it suffices to verify the first equality in \eqref{191215e}, or equivalently, that for every $a'\in M(x_1,x_2)$ and 
$b'\in M(x_2,x_3)$, there exists $b\in M(x_2,x_3)$ such that
\qtnl{081215b}
a\cdot b=a'\cdot b'.
\eqtn 
To do this, we observe that by the saturation property of $\cM$,
there exists a vertex $y\in X$ adjacent to each of the vertices $x_1$, 
$x_2$, $x_3$ in the graph~$D$. It follows 
that $x_1\sim y\sim x_2\sim x_1$ and $x_2\sim y\sim x_3\sim x_2$.
By formula~\eqref{230116a}, this implies that
\qtnl{051215s}
M(x_1,y)\cdot M(y,x_2)=M(x_1,x_2),\qquad M(x_2,y)\cdot M(y,x_3)=M(x_2,x_3).
\eqtn
Using these equalities, we successively find $u\in M(x_1,y)$ and 
$t\in M(x_2,y)$  such that $a'\cdot t=u$, and then  $v\in M(y,x_3)$ 
such that $t\cdot v=b'$. Then
\qtnl{241115w}
a'\cdot b'=(u\cdot t^*)\cdot(t\cdot v)=u\cdot v.
\eqtn
Using equalities~\eqref{051215s} again, we first find $s\in M(x_2,y)$ 
such that $a\cdot s=u$, and then $b\in M(x_2,x_3)$ such that $s^*\cdot b=v$
(the obtained configuration is depicted at Fig.~\ref{f5}). 
\begin{figure}[h]
\grphp{
& & & & & \VRT{x_2}\ar[dd]^{t}\ar[ddrrrrr]^{b'} & & & & &\\
& & & & & & & & & &\\
\VRT{x_1} \ar[uurrrrr]^{a'}\ar[rrrrr]_{u}\ar[ddrrrrr]_{a} & & & & & \VRT{y}\ar[rrrrr]_{v} & & & & & \VRT{x_3} \\
& & & & & & & &\\
& & & & & \VRT{x_2}\ar[uu]_{s}\ar[uurrrrr]_{b} & & & & &  \\
}
\caption{}\label{f5}
\end{figure} 
Thus, from~\eqref{241115w}, it follows that
$$
a'\cdot b'=u\cdot v=(a\cdot s)\cdot(s^*\cdot b)=a\cdot b,
$$
which proves~\eqref{081215b}. This completes the proof of~\eqref{191215e}, and, hence, the first statement.\medskip

To prove the second statement, let $P$ and $P'$ be a $d$- and a $d'$-path connecting the elements $u=x^{}_1=x'_1$ and $v=x^{}_{d^{}+1}=x'_{d'+1}$.
Without loss of generality, we can assume that $d+d'\ge 3$. 
By formula~\eqref{230116a}, the required statement follows for
$d+d'=3$. Suppose that $d+d'\ge 4$. Then since the matching configuration
$\cM$ is saturated, there exists  a vertex $w\in X$ adjacent with 
each of the vertices $\{u,x_2,x_3,x'_2\}$ in the graph~$D$.
By formula~\eqref{191215e} applied  
to $(u,w,x_2)$ and $(x_2,w,x_3)$, we have
$$
M_P=M(u,x_2)\cdot M(x_2,x_3)\cdot M_{Q_1}=(M(u,w)\cdot a)\cdot (a^*\cdot M(w,x_3))\cdot M_{Q_1}=
$$
$$
M(u,w)\cdot M(w,x_3)\cdot M_{Q_1}=M(u,w)\cdot M_Q,
$$
where $Q_1=(x_3,\ldots,v)$, $Q=(w,x_3,\ldots,v)$, and $a\in M(w,x_2)$.
Similarly, one can prove that $M_{P'}=M(u,w)\cdot M_{Q'}$, where
$Q'=(w,x'_2,\ldots,v)$. Note that $Q$ and $Q'$ are respectively, 
$(d-1)$- and $d'$-paths connecting $w$ and $v$. By the induction 
hypothesis, this implies that $M_{Q^{}}=M_{Q'}$ and hence
$$
M_P=M(u,w)\cdot M_Q=M(u,w)\cdot M_{Q'}=M_{P'}
$$
as required.\eprf\medskip

In what follows, wee keep the notation of Theorem~\ref{200116a}.
Let us represent the set $T'=(M\cdot M)^\#$ as the union of the sets
$$
T'(x,y)=\bigcup_{z\in N(x,y)}M(x,z)\cdot M(z,y),\quad x,y\in X,
$$
where $N(x,y)=\{z\in X:\ x\sim z\sim y\}$.
Note that the latter set is not empty for all $x$ and $y$, because the 
matching configuration~$\cM$ 
is saturated. Moreover, if $z\in N(x,y)$, then $M(x,z)\cdot M(z,y)=M_P$,
where $P=(x,z,y)$ is a $2$-path. By Lemma~\ref{191215b}, this 
implies that
$T'(x,y)$ is the partition of $\Delta_x\times\Delta_y$ into matchings.
In particular,
\qtnl{260216e}
1_x=u\cdot u^*\in M(x,z)\cdot M(z,x)=T'(x,x),
\eqtn
where $u\in M(x,z)$. 

\lmml{230216a}
The pair $\cY'=(\Delta,T')$ is a coherent configuration.
\elmm
\proof By the above remarks, the set $T'$ is a partition 
of $\Delta\times\Delta$ into matchings and $1_\Delta\in (T')^\cup$. 
It is easily seen that $(T')^*=T'$.  Thus, it suffices to verify that 
if $u,v\in T'$ and $u\cdot v\ne\varnothing$, then $u\cdot v\in T'$. 
To this end, let $u\in M_P$ and $v\in M_Q$, where $P$ and $Q$
are $2$-paths of the graph $D$. Since $u\cdot v\ne\varnothing$,
the last vertex of $P$ coincides with the first vertex of~$Q$.
By Lemma~\ref{191215b}, this implies that
$$
M_P\cdot M_Q=M_{P\cdot Q},
$$
where $P\cdot Q$ is the $4$-path of $D$ consisting of the vertices
of $P$ followed by the verices of $Q$ (the last vertex of $P$ 
is identified with the first vertex of~$Q$). Thus,
$u\cdot v\in M_{P\cdot Q}$ belongs to $T'$, as required.\eprf\medskip

To complete the proof of Theorem~\ref{200116a}, we note 
that the coherent configuration~$\cY'$ defined in Lemma~\ref{230216a}
contains the coherent closure $\cY$ of~$M$. Indeed, let $x\sim y$. Then 
$$
M(x,y)=1_x\cdot M(x,y)\subset T'.
$$
This implies that $M\subseteq T'\subseteq (T')^\cup$. Thus, the
claim follows, because~$\cY$ is the smallest coherent configuration,
for which every relation of~$M$ is the union of basis relations.
Conversely, as is easily 
seen, every relation in the set $M\cdot M$ is contained in $T^\cup$. 
By the definition of $T'$, this implies that $T'\subseteq T^\cup$, 
i.e., $\cY'\le\cY$. Thus, $\cY'=\cY$, as required.\eprf\medskip

The set $T$ in  Theorem~\ref{200116a} consists of matchings.
Therefore, $n_t=1$ for all $t\in T$. By 
formula~\eqref{260216e}, this proves the following statement.

\crllrl{260216ts}
In the notation of Theorem~\ref{200116a}, the coherent
configuration~$\cY$ is semiregular and $\Phi(\cY)=\Phi(\cM)$.
\ecrllr

From Corollary~\ref{260216ts} and Theorem~\ref{150316c}, it follows
that every saturated matching configuration as well as its coherent
closure are schurian. Since the latter one is also separable, the
following statement says that, in fact, every saturated matching
configuration is also separable (cf. Lemma~\ref{150316d}).

\thrml{260216t}
Let $\cM=(\Delta,M)$ and $\cM'=(\Delta',M')$ be saturated
matching configurations  with coherent closures $\cY=(\Delta,T)$ 
and $\cY'=(\Delta',T')$, respectively. Let 
$\varphi:M\to M',\ u\mapsto u'$ be 
an algebraic isomorphism from $\cM$ onto $\cM'$.  
Then 
$$
T=(M\cdot M)^\#\qaq T'=(M'\cdot M')^\#,
$$ 
and the mapping
\qtnl{280216a}
\psi:T\to T',\quad u\cdot v\mapsto u'\cdot v'
\eqtn
is a well-defined bijection. Moreover, $\psi|_M=\varphi$
and $\psi$ is an algebraic isomorphism from $\cY$ to $\cY'$.
\ethrm
\proof The first statement immediately follows from Theorem~\ref{200116a}. To verify that the mapping $\psi$ is
well-defined, suppose that $b_1\cdot c_1=b_2\cdot c_2$ 
for some $b_1,c_1,b_2,c_2\in M$. Let $x,y,z_1,z_2\in X$ be such that
$$
b_i\in M(x,z_i)\qaq c_i\in M(z_i,y),\qquad i=1,2.
$$
In particular, $x\sim z_i\sim y$ for each~$i$. Since the matching
configuration $\cM$ is saturated, there exists a vertex $z\in X$
adjacent to each of the vertices  $x,z_1,z_2,y$ in the graph~$D$.
By formula~\eqref{230116a}, this implies that 
\qtnl{290309g}
M(x,z_i)=M(x,z)\cdot M(z,z_i)\qaq M(z,y)=M(z,z_i)\cdot M(z_i,y)
\eqtn
for each~$i$. Take any $a_1\in M(x,z)$. Then the first equality in~\eqref{290309g} implies that 
\qtnl{260216c}
d_1:=a_1^*\cdot b_1\in M(z,z_1)\qaq d_2:=a_1^*\cdot b_2\in M(z,z_2),
\eqtn
whereas by the second equality in~\eqref{290309g}, we have
\qtnl{260216d}
a_2:=d_1\cdot c_1\in M(z,y).
\eqtn
Thus,
\qtnl{260216b}
a_1\cdot a_2=(a_1\cdot d_1)\cdot (d^*_1\cdot a_2)=b_1\cdot c_1=b_2\cdot c_2=a_1\cdot d_2\cdot c_2,
\eqtn
whence $c_2=d_2^*\cdot a_2$ (see Fig.~\ref{f59}).\medskip
\begin{figure}[h]
	\grphp{
		& & & \VRT{z_1} \ar[dddrrr]^*{c_1} & & & \\
		& & & & & & \\
		& & & & & & \\
		\VRT{x}\ar[rrruuu]^*{b_1}\ar[rrr]_*{a_1}\ar[rrrddd]_*{b_2} & & &
		\VRT{z}\ar[uuu]_*{d_1}\ar[rrr]_*{a_2}\ar[ddd]^*{d_2} & & &
		\VRT{y} \\
		& & & & & & \\
		& & & & & &\\
		& & & \VRT{z_2} \ar[uuurrr]_*{c_2} & & & \\
	}
	\caption{}\label{f59}
\end{figure}

On the other hand, the algebraic isomorphism $\varphi$ induces
a bijection from $\Phi(\cM)$ onto $\Phi(\cM')$ and hence a bijection
$x\mapsto x'$ from $X$ onto $X'$. In these notation,  
$$
M(x,y)^\varphi=M'(x',y').
$$ 
In particular, $x\sim y$ if and only
if $x'\sim y'$ for all $x,y\in X$. Thus, the bijection $x\mapsto x'$
is a graph isomorphism from $D(\cM)$ onto $D(\cM')$.
Moreover, equality~\eqref{230116a} implies that if
$x\sim y\sim z\sim x$ and $a\in M(x,y)$, $b\in M(y,z)$, 
then $(a\cdot b)'=a'\cdot b'$. Thus, in view of
formulas~\eqref{260216c}, \eqref{260216d}, and ~\eqref{260216b},
we have
$$
b'_1\cdot c'_1=(a_1\cdot d_1)'\cdot (d^*_1\cdot a_2)'=
a'_1\cdot (d'_1\cdot (d'_1)^*)\cdot a'_2=a'_1\cdot a'_2=
$$
$$
(b_2\cdot d_2^*)'\cdot (d_2\cdot c_2)'=
b'_2\cdot (d'_2\cdot (d^*_2)')\cdot c'_2=b'_2\cdot c'_2,
$$
which proves that $\psi$ is a well-defined bijection. Note that if 
$a\in  M(x,y)$ for some $x$ and $y$, then
$$
a^\psi=(1_x\cdot a)^\psi=1_{x'}\cdot a'=a'=a^\varphi.
$$
Thus, $\psi|_M=\varphi$.\medskip

To prove that $\psi$ is an algebraic isomorphism from $\cY$
onto $\cY'$, we note that if $a,b\in T$, then either
$a\cdot b\in T$ or $a\cdot b=\varnothing$. Thus,
it suffices to verify that
\qtnl{300309a}
(a\cdot b)^\psi=a^\psi\cdot b^\psi,\qquad a,b\in T,\ 
a\cdot b\ne\varnothing.
\eqtn
To this end, take such relations $a$ and $b$. Then 
$$
a\in M(x,z_1)\cdot M(z_1,z')\qaq b\in M(z',z_2)\cdot M(z_2,y)
$$
for appropriate $x,z_1,z',z_2,y\in X$. Since $\cM$ is saturated, we 
may assume that $z_1=z_2$; denote this element by $z$. Then 
$x\sim z'\sim z\sim x$ and $y\sim z'\sim z\sim y$. In view of
formula~\eqref{230116a}, one can find $c_1\in M(x,z)$, $c_2\in M(z,z')$, 
and  $c_3\in M(z,y)$ such that
$$
c_1\cdot c_2=a\qaq c_2\cdot b=c_3
$$
(see Fig.~\ref{f6}). 
\begin{figure}[h]
	\grphp{
		& & & \VRT{z} \ar[ddddrrr]^*{c_3}\ar[dddd]^*{c_2} & & & \\
		& & & & & & \\
		& & & & & & \\
		& & & & & & \\
		\VRT{x}\ar[rrruuuu]^*{c_1}\ar[rrr]_*{a} & & &
		\VRT{z'}\ar[rrr]_*{b} & & &
		\VRT{y} \\
	}
	\caption{}\label{f6}
\end{figure}
Besides, by the definition of $\psi$, we have 
$(c_1\cdot c_3)^\psi=c_1^\varphi\cdot c_3^\varphi$. 
Since $\psi|_M=\varphi$, we obtain
$$
(a\cdot b)^\psi=(c_1^{}\cdot c_2^{}\cdot c_2^*\cdot c_3^{})^\psi=
(c_1^{}\cdot c^{}_3)^\psi=
c_1^\varphi\cdot c_3^\varphi=c_1^\psi\cdot c_3^\psi=
$$
$$
c_1^\psi \cdot c_2^\psi\cdot (c_2^*)^\psi\cdot c_3^\psi=
(c_1^{}\cdot c_2^{})^\psi \cdot (c_2^*\cdot c_3^{})^\psi=
a^\psi\cdot b^\psi,
$$
which completes the proof of \eqref{300309a}.\eprf

\section{Matching subconfiguration of a coherent configuration}\label{220216c}

\sbsnt{Main statement.}\label{110316y}
Let $\cX=(\Omega,S)$ be a homogeneous coherent configuration with maximum
valency~$k$. Define a binary relation $\sim$ on the set $X=S_k$ by setting 
\qtnl{290216a}
x\sim y\qquad \Leftrightarrow\qquad |x^*y|=k. 
\eqtn 
Formally, this relation is different from that introduced in 
Section~\ref{030416a}, but as we will see, they are similar in some sense. 
The relation~\eqref{290216a} is symmetric, but, in general, not reflexive  
or transitive. Formula~\eqref{150410a} implies that 
$$
k^2=n_{x^*}n_{y^{}}=\sum_{z\in x^*y}n_zc_{x^*y^{}}^z=
\sum_{z\in x^*y}n_{y^*}c_{z^*x^*}^{y^*}=k\sum_{z\in x^*y}c_{xz}^y
$$
(here we also use that $n_{x^{}}=n_{x^*}=n_{y^{}}=n_{y^*}=k$).
Since $c_{xz}^y\ge 1$ for all $z\in x^*y$, we conclude that
\qtnl{290216d}
x\sim y\qquad \Leftrightarrow\qquad 
c_{xz}^y=1\ \,\text{for all}\ \, z\in x^*y. 
\eqtn

Fix a point $\alpha\in\Omega$. Then this formula implies also
that $x\sim y$ if and only if each nonempty relation $r_{x,y}$
defined in Lemma~\ref{260511a} is a matching. Set
\qtnl{240216a}
\Delta_\alpha=\alpha X\qaq 
M_\alpha=\bigcup_{x,y\in X} M_\alpha(x,y),
\eqtn
where 
\qtnl{260216v}
M_\alpha(x,y)=
\css
S_\alpha(x,y)      &\text{if $x\sim y$,}\\
\varnothing        &\text{if $x\not\sim y$ and $x\ne y$,}\\
\{1_{\alpha x}\}   &\text{if $x\not\sim y$ and $x=y$,}\\
\ecss
\eqtn
with $S_\alpha(x,y)=\{r_{x,y}:\ r\in x^*y\}$. Thus, 
the pair $\cM_\alpha=(\Delta_\alpha,M_\alpha)$ is close to
be a matching configuration, see Definition~\ref{290216b}.
The problem is that, in general, $\cM_\alpha$ is not 
necessarily a partial coherent configuration:
formula~\eqref{230116a} is not always true. 

\thrml{240116a}
Let $\cX=(\Omega,S)$ be a homogeneous coherent 
configuration with maximum valency~$k>1$ and  indistinguishing
number~$c$. Suppose that
\qtnl{240116b}
|X|>6c(k-1),
\eqtn
where $X=S_k$. Then given a point $\alpha\in\Omega$, the pair
$\cM_\alpha$ is a saturated matching configuration; in particular, 
$\Phi(\cM_\alpha)=\{\alpha x:\ x\in X\}$. Moreover,
\qtnl{050316c}
S_\alpha(x,y)\subseteq (M_\alpha(x,z)\cdot M_\alpha(z,y))^\cup 
\eqtn
for all $x,y,z\in X$ such that $x\sim z\sim y$. 
\ethrm

We prove Theorem~\ref{240116a} in Subsection~\ref{100416a}.
The key points of the proof are the lemmas  proved in
Subsections~\ref{240116f} and~\ref{130316a}.

\sbsnt{Neighborhoods.}\label{240116f}
In this subsection, we are interested in estimating the cardinality 
of the {\it neighborhood} $N(Y)$ of a set $Y\subseteq X$ with respect to the
relation~$\sim$ that is defined as follows:
\qtnl{250216a}
N(Y)=\{x\in X:\ x\sim y\ \,\text{for all}\ \,y\in Y\}.
\eqtn
In the sequel, we set $N(x,y,\ldots)=N(\{x,y,\ldots\})$. The following
statement gives a lower bound for the cardinality of the neighborhood.

\lmml{291115c}
If $Y\subseteq X$, then $|N(Y)|\ge|X|-c(k-1)|Y|$.
\elmm
\proof  Obviously,
$$
|N(Y)|\ge|X|-|Y|\,\max_{y\in Y}|X_y|, 
$$
where $X_y$ is the set of all $x\in X$ such that $y\not\sim x$.
Therefore, it suffices to verify that for any $y\in X$,
\qtnl{240616a}
|X_y|\le c(k-1).
\eqtn
To do this fix a relation $y\in X$, a point $\alpha\in \Omega$,
and denote by $\Lambda_y$ the set of all pairs of distinct points
of~$\alpha y$. From~\eqref{290216d},  it follows that 
$x\in X_y$ only if $c_{yz}^x>1$ for some $z\in y^*x$. 
In the latter case, for each $\beta\in\alpha x$, there exists
a pair $(\gamma,\delta)\in\Lambda_y$ such that 
\qtnl{050316a}
r(\gamma,\beta)=z=r(\delta,\beta).
\eqtn
It follows that the set $T_{x,y}$ of all such triples
$(\beta,\gamma,\delta)$ contains at least $|\alpha x|=n_x=k$ elements.
Therefore, the union of all sets $T_{x,y}$ with $x\in X_y$ contains
at least $k|X_y|$ elements. On the other hand,
$$
\bigcup_{x\in X_y}T_{x,y}=
\bigcup_{e\in \Lambda_y}T_e,
$$
where $T_e$ is the set of all triples belonging to the left-hand side
and the second and third entries of which forms the pair~$e$. Moreover, 
the number of nonempty summands on the right-hand side is at most
$|\Lambda_y|=k(k-1)$. Thus, there exists a pair $e\in\Lambda_y$ such that 
$$
|T_e|\ge \frac{1}{k(k-1)}\sum_{x\in N_y}|T_{x,y}|\ge
\frac{k|X_y|}{k(k-1)}.
$$
Note that in  view of~\eqref{050316a} if $e=(\gamma,\delta)$, then  
$T_e$ is contained in the set $\Omega_{\gamma,\delta}$ defined in
formula~\eqref{100814c}. Consequently,
$$
c\ge|\Omega_{\gamma,\delta}|\ge |T_e|\ge\frac{|X_y|}{k-1}.
$$
which proves formula~\eqref{240616a}.\eprf

\sbsnt{Special elements.}\label{130316a}
Let $x,y,z\in X$ and $r,s,t\in S$  be such that 
\qtnl{240616t}
x\sim z\sim y\qaq r\in x^*z,\ \,s\in z^*y,\ \,t\in x^*y.
\eqtn
Then the relations $r_{x,z}\in S_\alpha(x,z)$ and 
$s_{z,y}\in S_\alpha(z,y)$ are matchings, whereas the relation
$t_{x,y}\in S_\alpha(x,y)$ is not necessarily a matching. 

\dfntn
An element $q\in N(x,y,z)$ is said to be {\it special} with respect 
to the $6$-tuple $T=(x,r,z,s,y,t)$, or {\it $T$-special} if 
there exist elements $u\in x^*q$, $v\in  z^*q$, and $w\in y^*q$ such that 
\qtnl{120316b}
uv^*\cap x^*z=\{r\}\qaq
vw^*\cap z^*y=\{s\}\qaq
uw^*\cap x^*y=\{t\},
\eqtn
see the configuration depicted in Fig.\ref{f667}.
\edfntn
\begin{figure}[h]
	\grphp{
		& & & \VRT{z} \ar[ddd]^*{v}\ar[ddddddrrr]^*{s} & & & \\
		& & & & & & \\
		& & & & & & \\
		& & & \VRT{q} & & & \\
		& & & & & & \\
		& & & & & &\\
		\VRT{x} \ar[uuuuuurrr]^*{r}\ar[rrruuu]_*{u}\ar[rrrrrr]_*{t}
		& & &  & & &
		\VRT{y} \ar[llluuu]^*{w}\\
	}
	\caption{}\label{f667}
\end{figure}

It follows from the definition that $q\sim x$, $q\sim y$, and $q\sim z$.
The following two statements reveal the property of special elements, 
to be used in the sequel, and provide a sufficient condition for their
existence.

\lmml{090316u}
In the above notation, assume that the set $N(x,y,z)$ contains 
a $T$-special element. Then $r_{x,z}\cdot s_{z,y}\subseteq t_{x,y}$.
\elmm
\proof Let $q$ be a $T$-special element. Then there exist elements 
$u\in x^*q$, $v\in  z^*q$, and $w\in y^*q$, for which
relations~\eqref{120316b} hold. It follows that
\qtnl{120316d}
u^{}_{x,q}\cdot v^*_{q,z}\subseteq r^{}_{x,z},\quad
v^*_{z,q}\cdot w^*_{q,y}\subseteq s^{}_{z,y},\quad
u^{}_{x,q}\cdot  w^*_{q,y}\subseteq t^{}_{x,y}.
\eqtn
On the other hand, the relations $u^{}_{x,q}\cdot v^*_{q,z}$ 
and $r^{}_{x,z}$ are matchings, because $x\sim q\sim z$, and 
$x\sim z$. Therefore, by
the first inclusion in~\eqref{120316d}, we conclude that
$u^{}_{x,q}\cdot v^*_{q,z}=r^{}_{x,z}$. Similarly,
$v^*_{z,q}\cdot w^*_{q,y}=s^{}_{z,y}$. Thus, by the third inclusion
in~\eqref{120316d}, we have
$$
r^{}_{x,z}\cdot s^{}_{z,y}=(u^{}_{x,q}\cdot v^*_{q,z})\cdot
(v^*_{z,q}\cdot w^*_{q,y})=u^{}_{x,q}\cdot  w^*_{q,y}
\subseteq t^{}_{x,y},
$$
which proves the required inclusion.\eprf

\lmml{110316x}
Let $T=(x,r,z,s,y,t)$, where $x,y,z\in X$ and $r,s,t\in S$ are such that
formula~\eqref{240616t} holds. Assume that 
\qtnl{240616q}
|N(x,y,z)|>3c(k-1)\qaq
(r_{x,z}\cdot s_{z,y})\,\cap\, t_{x,y}\ne\varnothing.
\eqtn
Then the set $N(x,y,z)$ contains a  $T$-special element. 
\elmm
\proof By the right-hand side relation in~\eqref{240616q},
there exist points $\beta\in \alpha x$, $\gamma\in\alpha z$,
and $\delta\in\alpha y$ such that
\qtnl{240616x}
(\beta,\gamma)\in r,\qquad (\gamma,\delta)\in s,\qquad
(\beta,\delta)\in t.
\eqtn
For every relation $q\in N(x,y,z)$ and every point $\mu\in\alpha q$
set
\qtnl{120316t}
u(q,\mu):=r(\beta,\mu),\qquad
v(q,\mu):=r(\gamma,\mu),\qquad
w(q,\mu):=r(\delta,\mu).
\eqtn
Then, obviously,
\qtnl{240616y}
r\in u v^*\, \cap\,  x^*z,\quad
s\in v w^*\, \cap\,  z^*y,\quad
t\in u w^*\, \cap\,  x^*y.
\eqtn
where $u=u(q,\mu)$, $v=v(q,\mu)$, and $w=w(q,\mu)$. Therefore,
if $q$ is not $T$-special, then for each $\mu\in \alpha q$,
there exists a basis relation $a=a(q,\mu)$ such that
\qtnl{240616z}
a\in (uv^*\cap x^*z)\setminus\{r\}\qoq
a\in (vw^*\cap z^*y)\setminus\{s\}\qoq
a\in (uw^*\cap x^*y)\setminus\{t\}.
\eqtn
Note that each of the sets $x^*z$, $z^*y$, $x^*y$ consists of at 
most $k$ relations. So  the relation $a(q,\mu)$
is one of the $3(k-1)$ relations contained in the
set $x^*z\cup z^*y\cup x^*y$, which does not depend on $q$ and $\mu$.
Thus, if $S_T$ is the set of all non-$T$-special relations~$q$ and
$P_a$ is the set of all pairs $(q,\mu)\in S_T\times \alpha q$
with $a=a(q,\mu)$, then there exists $a\in x^*z\cup z^*y\cup x^*y$ such that
\qtnl{120316w}
|P_a|\ge \frac{k|S_T|}{3(k-1)}.
\eqtn
Without loss of generality, we may assume that
$a\in x^*z$. Let now $(q,\mu)\in P_a$. Then by the definition
of $u=u(q,\mu)$ and $v=v(q,\mu)$, we have 
$\mu\in\beta u\cap\gamma v$. Since also $a\in uv^*$ (see~\eqref{240616z}),
there exists a point $\nu(q,\mu)$ belonging to the set 
$\mu u^*\cap \gamma a^*$,
see the configuration depicted in Fig.~\ref{f158}. 
\begin{figure}[h]
	\grphp{
		& & & & \VRT{\mu} & &  & &
		\VRT{\nu}\ar[llll]_*{u} \ar[llllddd]^*{a} \\
		& &  & & & &  & &\\
		& &  & & & &  & &\\
		\VRT{\beta}\ar[rrrruuu]^*{u}\ar[rrrr]_*{r} & & & &	
		\VRT{\gamma}\ar[uuu]_*{v} & & & &\\		
	}
	\caption{}\label{f158}
\end{figure} 
Note that every point $\nu(q,\mu)$ belongs to the set $\gamma a^*$
of cardinality at most~$k$. Therefore, there exists a point 
$\nu\in \gamma a$ such that the set $P_{a,\nu}$ of all pairs 
$(q,\mu)\in P_a$ with $\nu=\nu(q,\mu)$ contains at least 
$|P_a|/k$ elements. Taking into account that $\mu$ is contained
in the set $\Omega_{\beta,\nu}$  defined by~\eqref{100814c},
we conclude by~\eqref{120316w} that
$$
c\ge|\Omega_{\beta,\nu}|=|P_{a,\nu}|\ge
\frac{|P_a|}{k}\ge \frac{k|S_T|}{3k(k-1)},
$$
whence $|S_T|\le  3c(k-1)$. By the hypothesis of the lemma, this implies
that the set $N(x,y,z)\setminus S_T$ is not empty. Since every element
of this set is $T$-special, we are done.\eprf

\sbsnt{Proof of Theorem~\ref{240116a}.}\label{100416a}
Let $x,y,z\in X$ be such that $x\sim y\sim z$. Then given 
matchings $a\in M(x,z)$ and $b\in M(z,y)$, there exists a
relation $c\in S(x,y)$ such that $a\cdot b$ intersects~$c$. It
follows that 
$$
(r_{x,z}\cdot s_{z,y})\cap t_{x,y}\ne\varnothing,
$$
where the relations $r,s,t\in S$ are defined by the conditions
$a=r_{x,z}$, $b=s_{z,y}$, and $c=t_{x,y}$, respectively. On the
other hand by the theorem hypothesis and Lemma~\ref{291115c},  
$$
|N(x,y,z)|\ge |X|-3c(k-1)>3c(k-1).
$$
Thus, the hypothesis of Lemma~\ref{110316x} holds for 
$T=(x,r,z,s,y,t)$. This shows that the set $N(x,y,z)$ contains a 
$T$-special element. By Lemma~\ref{090316u}, this implies that
$$
a\cdot b=r_{x,z}\cdot s_{z,y}\subseteq t_{x,y}.
$$
This proves formula~\eqref{050316c}. Moreover, if $x\sim z$, then
$t_{x,y}$ is a matching of size~$k$. Since $a\cdot b$ is also  a
matching of the same size, we conclude that $a\cdot b=c$. Thus,
condition~\eqref{230116a} is satisfied and hence $\cM_\alpha$ is 
a matching configuration. To prove that it is saturated, it
suffices to verify that the set $N(Y)$ is not empty for all 
$Y\subseteq X$ with $|Y|\le 4$. But in this case,
$$
|N(Y)|\ge |X|-4c(k-1)>2c(k-1)>0
$$
by formula~\eqref{240116b} and Lemma~\ref{291115c}.\eprf

\sbsnt{Algebraic isomorphisms.}
In this subsection, we fix an algebraic isomorphism 
$\varphi:\cX\to\cX'$, $r\mapsto r'$,
where $\cX$ is a coherent configuration satisfying
condition~\eqref{240116b} and $\cX'=(\Omega',S')$ is an arbitrary
coherent configuration.

\thrml{140316a}
In the above notation and assumptions, the coherent configuration~$\cX'$
satisfies condition~\eqref{240116b}. Moreover, given points 
$\alpha\in \Omega$ and $\alpha'\in\Omega'$, the mapping 
$$
\varphi_{\alpha,\alpha'}:M^{}_{\alpha^{}}\to M'_{\alpha'},\quad
r^{}_{x^{},y^{}}\mapsto r'_{x',y'}
$$
is an algebraic isomorphism from $\cM_\alpha$ to $\cM'_{\alpha'}$.
Moreover, if $a\in M_\alpha$ and $a\subset r\in S$, then 
$a^{\varphi_{\alpha,\alpha'}}\subset r^\varphi$.
\ethrm
\proof Clearly, the rank, maximal valency, and indistinguishing number
of $\cX'$ are the same as of $\cX$. By Theorem~\ref{240116a}, given
$\alpha'\in\Omega'$, one can define a saturated
matching configuration $\cM'_{\alpha'}=(\Delta'_{\alpha'},M'_{\alpha'})$
exactly in same way as for $\cX$. Since the algebraic isomorphism
$\varphi$ takes $X:=S^{}_k$ to $X':=S'_k$ and $r_{x,y}$ is not empty 
if and only if so is $r'_{x',y'}$, the mapping
$\psi=\varphi_{\alpha,\alpha'}$ is a bijection. We
need the following lemma.

\lmml{070316a}
Let  $x,y,z\in X$ be such that $x\sim y\sim z$, and let 
$r\in x^*z$, $s\in z^*y$. Then $r_{x,z}\cdot s_{z,y}\subseteq t_{x,y}$
for some $t\in x^*y$. Moreover
$$
(r^{}_{x^{},z^{}})^\psi\cdot (s^{}_{z^{},y^{}})^\psi
\subseteq t'_{x',y'}.
$$
\elmm
\proof Let $a=r_{x,z}$, $b=s_{z,y}$, and $c=a\cdot b$. Then
$c\in M_\alpha(x,z)\cdot M_\alpha(z,y)$. By formula~\eqref{050316c},
this implies that $c\subseteq t_{x,y}$, where $t$ is a unique 
relation in $x^*y$ that 
intersects~$c$. This proves the first statement. To prove the second one,
we note that 
$$
r_{x,z}\cdot s_{z,y}\,\cap\, t_{x,y}\ne\varnothing
\qaq|N(x,y,z)|>3c(k-1),
$$
where the latter follows from 
Lemma~\ref{291115c}. In view of Lemma~\ref{110316x}, this implies
that $N(x,y,z)$ contains a  $T$-special element $q$, where
$T=(x,r,z,s,y,t)$. This means that formula~\eqref{120316b} holds
and hence
$$
u'(v')^*\cap (x')^*z'=\{r'\}\qaq
v'(w')^*\cap (z')^*y'=\{s'\}\qaq
u'(w')^*\cap (x')^*y'=\{t'\}
$$
for suitable relations $u'\in (x')^*q'$, $v'\in (z')^*q'$, and 
$w'\in (y')^*q'$. Therefore, $q'$ is a $T'$-special element, where
$T'=(x',r',z',s',y',t')$. By Lemma~\ref{090316u}, this shows
that  
$$
r'_{x',z'}\cdot s'_{z',y'}\subseteq t'_{x',y'}, 
$$
as required.\eprf\medskip

To complete the proof of Theorem~\ref{140316a}, assume that in Lemma~\ref{070316a}, $x\sim y$.
Then $r_{x,z}\cdot s_{z,y}=t_{x,y}$ is a matching. Since
$(t_{x,y})'=t'_{x',y'}$ by the definition of $\varphi_{\alpha,\alpha'}$,
we obtain
$$
(r^{}_{x^{},z^{}})^\psi\cdot (s^{}_{z^{},y^{}})^\psi=
(r^{}_{x^{},z^{}}\cdot s^{}_{z^{},y^{}})^\psi.
$$
Thus, $\psi$ is an algebraic isomorphism from $\cM^{}_{\alpha^{}}$
onto $\cM'_{\alpha'}$, as required.\eprf

\section{Proof of Theorem~\ref{170116c}}\label{090316r}
By the theorem hypothesis, the set $X:=S_k$ coincides with 
$S\setminus S_1$. Since also $|S|$ is greater than $m+6c(k-1)$, we have
\qtnl{260216u}
|X|=|S_k|=|S|-|S_1|>m+6c(k-1)-m=6c(k-1),
\eqtn
i.e., inequality~\eqref{240116b} holds. By Theorem~\ref{240116a},
this implies that for each point $\alpha\in\Omega$, the pair
$\cM_\alpha=(\Delta_\alpha,M_\alpha)$ defined by~\eqref{240216a} 
is a saturated matching configuration with fibers $\Omega_x=\alpha x$,
$x\in X$, and formula~\eqref{050316c} holds. 

\lmml{260216r}
Let $\cY_\alpha=(\Delta_\alpha,T_\alpha)$ be the coherent closure 
of $M_\alpha$, and let $\cD_\alpha$ be the 
complete coherent configuration on $\Omega\setminus\Delta_\alpha$. Then
\qtnl{250216y}
\cX_\alpha=\cD_\alpha\,\boxplus\, \cY_\alpha.
\eqtn
In particular, $\Phi(\cX_\alpha)=\{\alpha s:\ s\in S\}$.
\elmm
\proof  By the definition of $M_\alpha$, we have $M_\alpha\subseteq (S_\alpha)^\cup$, where $S_\alpha$ is the set
of basis relations of the coherent configuration~$\cX_\alpha$. 
By formula~\eqref{050316c} and Theorem~\ref{200116a}, this implies that
$$
T_\alpha=(M_\alpha\cdot M_\alpha)^\#\subset M_\alpha\cdot M_\alpha\subseteq (S_\alpha)^\cup.
$$
Thus, $(\cX_\alpha)_{\Delta_\alpha}\ge \cY_\alpha$ by the minimality 
of the coherent configuration $\cY_\alpha$.  Moreover, 
$\Omega\setminus\Delta_\alpha=\alpha S_1$ and $\alpha s\in\Phi(\cX_\alpha)$
for all $s\in S_1$. Therefore 
$$
\cX_\alpha=
(\cX_\alpha)_{\Omega\setminus\Delta_\alpha}\,\boxplus\,
(\cX_\alpha)_{\Delta_\alpha}\ge
\cD_\alpha\,\boxplus\, \cY_\alpha\ge \cX,
$$
and we are done by the minimality of the coherent configuration
$\cX_\alpha$.\eprf

\crllrl{150316a}
The coherent configuration $\cX_\alpha$ is schurian and separable.
\ecrllr
\proof By Corollary~\ref{260216t}, the coherent configuration
$\cY_\alpha$ is semiregular. By \eqref{250216y}, this implies that
$\cX_\alpha$ is $1$-regular. Thus, we are done by 
Theorem~\ref{150316c}.\eprf\medskip

Let us prove that $\cX$ is separable. By Lemma~\ref{150316d} and
Corollary~\ref{150316a}, it suffices to verify that every 
algebraic isomorphism $\varphi:\cX\to\cX'$ can be extended to
an algebraic isomorphism
\qtnl{190316a}
\varphi_0:\cX^{}_{\alpha^{}}\to\cX'_{\alpha'}.
\eqtn 
To construct this extension, we note that by Theorem~\ref{140316a},
given $\alpha\in\Omega'$, the mapping 
$\varphi_{\alpha,\alpha'}:r_{x,y}\mapsto r'_{x',y'}$, is
an algebraic isomorphism between the saturated matching configurations
$\cM^{}_{\alpha^{}}$ and $\cM'_{\alpha'}$. By Theorem~\ref{260216t}, 
this algebraic isomorphism can be extended to the algebraic isomorphism
\qtnl{140316b}
\psi_{\alpha,\alpha'}:\cY^{}_{\alpha^{}}\to \cY'_{\alpha'},
\eqtn
where $\cY^{}_{\alpha^{}}$ and $\cY'_{\alpha'}$ are the coherent 
closures of $\cM^{}_{\alpha^{}}$ and $\cM'_{\alpha'}$, respectively.
It remains to note that by Lemma~\ref{260216r}, the
algebraic isomorphism~\eqref{140316b}  can be extended to the
algebraic isomorphism~\eqref{190316a}
by setting
$$
a^{\varphi_0}=\css
a^{\psi_{\alpha,\alpha'}}         &\text{if $a\in M_\alpha$,}\\
r^\varphi_{s^\varphi,t^\varphi}   &\text{if $a\not\in M_\alpha$
	and $a=r_{s,t}$,}\\
\ecss
$$
here in view of formula~\eqref{250216y}, every basis relation 
of $\cX_\alpha$ either belongs to $M_\alpha$ or is of the form 
$r_{s,t}$ with $s,t\in S$.\medskip
 
To prove that $\cX$ is schurian, fix $\alpha\in\Omega$ and set
$\varphi=\id_S$. By the result of the previous paragraph with
$\cX'=\cX$ and $\alpha'\in\Omega$, the algebraic isomorphism 
$\varphi$ can be extended to the algebraic 
isomorphism~\eqref{190316a}. Therefore the set
$$
G_{\alpha\mapsto\alpha'}=
\iso(\cX_{\alpha^{}},\cX_{\alpha'},\varphi_0)
$$
is not empty by Corollary~\ref{150316a}. On the other hand,
 formula~\eqref{190316b} shows this set is contained in
the group $\iso(\cX,\cX,\id_S)=\aut(\cX)$.
Thus, the latter group  is transitive, because it contains
a transitive subgroup generated by 
the sets $G_{\alpha\mapsto\alpha'}$, $\alpha'\in\Omega$.
Moreover, 
$$
\orb(\aut(\cX)_\alpha)=\{\alpha s:\ s\in S\}
$$
by the second part of Lemma~\ref{260216r} and Corollary~\ref{150316a}.
Therefore, the coherent configuration associated with $\aut(\cX)$ 
coincides with $\cX$, i.e., $\cX$ is schurian.

\end{document}